\documentclass{amsart}
\usepackage[english]{babel}
\usepackage[latin1]{inputenc}
\usepackage[dvips,final]{graphics}
\usepackage{amsmath,amsfonts,amssymb,amsthm,amscd,array,stmaryrd,mathrsfs}
\usepackage{pstricks}
 \usepackage[all]{xy}
\usepackage{textcomp}
 \usepackage[final]{epsfig}
\theoremstyle{plain}
\newtheorem{thm}{Theorem}
\newtheorem{lem}{Lemma}[section]
\newtheorem{cor}[lem]{Corollary}
\newtheorem{prop}[lem]{Proposition}

\theoremstyle{definition}
\newtheorem{defn}[lem]{Definition}
\newtheorem{rem}[lem]{Remark}
\newtheorem{ex}[lem]{Example}

\newcommand{\R}{\mathbb{R}}
\newcommand{\Z}{\mathbb{Z}}
\newcommand{\C}{\mathbb{C}}
\newcommand{\N}{\mathbb{N}}

\newcommand{\Kb}{\mathbb{K}}

\newcommand{\Af}{\textbf{A}}
\newcommand{\Jf}{\textbf{J}}
\newcommand{\Kf}{\textbf{K}}

\newcommand{\Rc}{\mathcal{R}}

\newcommand{\Qc}{\mathcal{R}}
\newcommand{\Pc}{\mathcal{P}}

\newcommand{\Sc}{\mathcal{S}}
\newcommand{\Dc}{\mathcal{D}}
\newcommand{\A}{\mathcal{A}}

\newcommand{\cc}{\Gamma}
\newcommand{\K}{\mathcal{K}}

\newcommand{\E}{\mathcal{E}}
\newcommand{\U}{\mathcal{U}}
\newcommand{\T}{\mathcal{T}}

\newcommand{\gog}{\mathfrak{g}}

\newcommand{\ga}{\mathfrak{a}}

\newcommand{\id}{\textup{Id}}

\newcommand{\Gr}{\textup{Gr}}

\newcommand{\Span}{\textup{Span}}

\newcommand{\im}{\textup{Im}}

\newcommand{\half}{\textstyle \frac{1}{2}}
\newcommand{\quart}{\textstyle \frac{1}{4}}
\newcommand{\End}{\textup{End}}

\newcommand{\Ker}{\textup{Ker}}
\newcommand{\e}{\varepsilon}
\newcommand{\wt}{\textup{wt}}
\newcommand{\asl}{\mathrm{asl}}
\newcommand{\osp}{\mathrm{osp}}
\newcommand{\Der}{\textup{Der}}
\newcommand{\ad}{\mathrm{ad}}

\def\a{\alpha}

\def\e{\varepsilon}

\def\r{\rho}
\def\s{\sigma}

\begin{document}

\title{Universal enveloping algebras of Lie antialgebras}

\author{S{\'e}verine Leidwanger
and  Sophie Morier-Genoud}

\address{
S{\'e}verine Leidwanger,
Institut Math\'e\-ma\-tiques de Jussieu,
Th\'eorie des groupes,
Universit\'e Denis Diderot Paris 7,
case 7012,
75205 Paris Cedex 13, France}

\email{leidwang@math.jussieu.fr}

\address{
Sophie Morier-Genoud,
Institut Math\'e\-ma\-tiques de Jussieu,
Universit\'e Pierre et Marie Curie Paris 6,
175 rue du Chevaleret, 75013 Paris, France}

\email{sophiemg@math.jussieu.fr}

\date{}


\begin{abstract}
Lie antialgebras is a class of supercommutative algebras recently appeared in symplectic geometry. 
We define the notion of enveloping algebra of a Lie antialgebra and study its properties.
We show that every Lie antialgebra is canonically related to a Lie superalgebra and
prove that its enveloping algebra is a quotient of the enveloping algebra of the corresponding Lie superalgebra.
\end{abstract}

\maketitle

\section*{Introduction}

In 2007, Ovsienko \cite{Ovsienko} introduced  a class of non-associative superalgebras called Lie antialgebras. The axioms were established after encountering two ``unusual'' algebraic structures in a context of symplectic geometry: the algebras, named $\asl_2$ and the conformal antialgebra $\A\K(1)$.
Ovsienko found that both of them are related to an odd bivector fields on $\mathbb{R}^{2|1}$ invariant under the action of the orthosymplectic Lie superalgebra $\osp(1|2)$.
The algebra $\A\K(1)$ is also related to the famous Neveu-Schwartz Lie superalgebra. 

It turns out that the Lie antialgebras $\asl_2$ and $\A\K(1)$ are particular cases of Jordan superalgebras, known as tiny Kaplansky superalgebra $K_3$ and full derivation superalgebra, respectively. 
The algebra $K_3$ appears in the classifications of Jordan superalgebras made in the middle of 70's by Kac \cite{Kac} and Kaplansky \cite{Kaplansky}, while the full derivation superalgebra is the main example of the article of McCrimmon on Kaplansky superalgebras \cite{McC}.
The axioms of Lie antialgebras appear throughout \cite{McC} but no general theory is developed. 
Lie antialgebras are somehow a combination of Kaplansky and Jordan superalgebras.

The basics for the representation theory of Lie antialgebras and their relation to Lie superalgebras was also developed in \cite{Ovsienko}. The theory is specific for Lie antialgebras and cannot be applied to the more general class of Jordan superalgebras.

The present paper is a representation theoretic approach to Lie antialgebras in continuation of~\cite{Ovsienko}.
We show that every Lie antialgebra canonically corresponds to a Lie superalgebra,
this construction was given in \cite{Ovsienko} without proof.
We define the notion of universal enveloping algebra of a Lie antialgebra and establish a universal property.  We show that the PBW property is not satisfied in general. 
Our main result is a realization of  the enveloping algebra of a Lie antialgebra as a quotient of the universal enveloping algebra of the corresponding Lie superalgebra. 
We deduce that every representation of a Lie antialgebra can be extended to a representation of the corresponding Lie superalgebra as announced in \cite{Ovsienko}.

Representations and modules of Jordan superalgebras, in particular of the tiny Kaplansky superalgebra, have been studied by several authors \cite{MZ1,MZ2,T1} (see also~\cite{MG}).
In this paper, we compare different approaches and definitions.

The paper is organized as follows.

In Section \ref{JKO}, we recall the definitions of Jordan superalgebras and Kaplansky superalgebras. We explain the relation between these algebras and the Ovsienko Lie antialgebras. The relation is based on results appearing in \cite{McC}. We also recall the notions of representations and modules for these algebras.

In Section \ref{UEA}, we define the universal enveloping algebra $\U(\ga)$ of a Lie antialgebra. We establish a universal property and discuss the Poincar\'e-Birkhoff-Witt theorem for $\U(\ga)$. The PBW property fails in general but holds true for $K_3$.

In Section \ref{liens}, we recall the construction \cite{Ovsienko} of the adjoint Lie superalgebra $\gog_\ga$ associated to~$\ga$. 
We formulate and prove our main theorem stating that the enveloping algebra $\U(\ga)$ is a quotient of $\U(\gog_\ga)$. As a consequence, we obtain the relation between representations  of $\ga$ and those of $\gog_\ga$ as announced in \cite{Ovsienko}.
We illustrate the constructions in the case of the tiny Kaplansky algebra.  

Section \ref{KAK1} is entirely devoted to the case of the algebra $\A\K(1)$, the conformal antialgebra, also known as full derivation algebra. We study its enveloping algebra and the representations.

Section \ref{ProofSect} contains the computations establishing the Jacobi identity in $\gog_\ga$. The computations are not straightforward and were missing in \cite{Ovsienko}.\\

\noindent
\textbf{Acknowledgments:} 
We are grateful to R. Berger, C. Conley and V.~Ovsienko for enlightening discussions and helpful references. We also would like to thank Marie Kreusch for useful comments.

\section{Jordan, Kaplansky and Ovsienko superalgebras}\label{JKO}

\subsection{Definitions}
First recall that a \textit{superalgebra} $\Af$ is an algebra with a grading over $\Z_2$, \textit{i.e}
$$
\Af=\Af_0\oplus \Af_1,  \quad \Af_i\cdot \Af_j\subset \Af_{i+j}.
$$
Elements in the subspaces $\Af_i$, $i=0,1$ are called homogeneous elements. For $a\in \Af_i$, the value $i$ is called the parity of $a$ that we denote by  $\bar{a}\in \{0,1\}$.

A superalgebra is called \textit{supercommutative} if it satisfies
$$
a\cdot b=(-1)^{\bar{a}\bar{b}}\,b\cdot a,
$$
for all homogeneous elements $a,b$.

A superalgebra is called \textit{half-unital} (cf \cite{McC}) if it contains an even element $\e$ such that
$$
\e\cdot a=a, \;a\in \Af_0, \quad \e\cdot a=\textstyle \half a, \;a \in \Af_1.
$$

Throughout this paper we consider algebras over the field $\Kb=\R$ or $\C$.\\

\noindent
\textbf{Lie antialgebras.} A \textit{Lie antialgebra} is a superalgebra $\ga=\ga_0\oplus\ga_1$
with a supercommutative product satisfying the following cubic identities:
\begin{enumerate}
\item[(LA0)] associativity of $\ga_0$
$$x_1\cdot (x_2 \cdot x_3)=(x_1\cdot x_2)\cdot x_3,$$
for all $x_1, x_2, x_3 \in \ga_0$,
\item[(LA1)] half-action
\begin{equation*}
x_1\cdot\left(x_2\cdot{}y\right)=
\half \left(x_1\cdot{}x_2\right)\cdot{}y,
\end{equation*}
for all $x_1,x_2\in\ga_0$ and $y\in\ga_1$, 


\item[(LA2)] Leibniz identity
\begin{equation*}\label{eq3}
x\cdot\left(y_1\cdot{}y_2\right)=
\left(x\cdot{}y_1\right)\cdot{}y_2+
y_1\cdot\left(x\cdot{}y_2\right),
\end{equation*}
for all $x\in\ga_0$ and $y_1,y_2\in\ga_1$, 

\item[(LA3)] odd Jacobi identity
\begin{equation*}\label{eq4}
y_1\cdot\left(y_2\cdot{}y_3\right)+
y_2\cdot\left(y_3\cdot{}y_1\right)+
y_3\cdot\left(y_1\cdot{}y_2\right)=0,
\end{equation*}
for all $y_1,y_2,y_3\in\ga_1$.

\end{enumerate}

A weak version of Lie antialgebras is also considered by replacing the axiom (LA1) by:

(LA1') ``full action''
$$x_1\cdot\left(x_2\cdot{}y\right)+x_2\cdot\left(x_1\cdot{}y\right)=\left(x_1\cdot{}x_2\right)\cdot{}y,$$
for all $x_1,x_2\in\ga_0$ and $y\in\ga_1$.\\


Note that, in the case where the even part of a Lie antialgebra is generated by the odd part, the associativity axiom (LA0), is follows from the axioms (LA1)--(LA3) (cf. Section \ref{assoc}).\\

\noindent
\textbf{Jordan (super)algebras.} 
An algebra $\Jf$ is a \textit{Jordan algebra} if it satisfies
\begin{enumerate}
\item[(J1)] commutativity, $a\cdot b=b\cdot a,$  for all $a,b \in \Jf$,
\item[(J2)] Jordan identity, $a^2\cdot (b\cdot a)=(a^2\cdot b)\cdot a, $  for all $a,b \in \Jf$.
\end{enumerate}
A superalgebra $\Jf=\Jf_0\oplus \Jf_1$ is a \textit{Jordan superalgebra} if it satisfies
\begin{enumerate}
\item[(SJ1)] supercommutativity,
\item[(SJ2)] super Jordan identity 
\begin{eqnarray*}
&&(a\cdot b)\cdot(c\cdot d)+(-1)^{\bar{b}\bar{c}}\,(a\cdot c)\cdot(b\cdot d)+ (-1)^{(\bar{b}+\bar{c})\bar{d}}(a\cdot d)\cdot(b\cdot c)\\
&=& ((a\cdot b)\cdot c)\cdot d+  (-1)^{(\bar{b}+\bar{c})\bar{d}+\bar{b}\bar{c}}((a\cdot d)\cdot c)\cdot b+(-1)^{(\bar{b}+\bar{c}+\bar{d})\bar{a}+\bar{c}\bar{d}}((b\cdot d)\cdot c)\cdot a,
\end{eqnarray*}
\noindent
for all $a,b,c,d $ homogeneous elements in $ \Jf$.\\
\end{enumerate}

We refer to \cite{Martinez} (see also the references therein) for the general theory of Jordan superalgebras.\\

\noindent
\textbf{Kaplansky superalgebras.} 
In an unpublished work \cite{Kaplansky}, Irving Kaplansky considered the following class of half-unital supercommutative superalgebras:
\begin{enumerate}
\item[\textbullet] let $\Kf_0$ be a unital commutative associative algebra,
\item[\textbullet] let $\star:  \Kf_0\times \Kf_1\rightarrow \Kf_1$ be a representation of $\Kf_0$,
\item[\textbullet] let $<\, , \, >: \Kf_1\times \Kf_1\rightarrow \Kf_0$ be  a $\Kf_0$-valued skew-symmetric bilinear form,
\end{enumerate}
define a structure of algebra on the space $\Kf=\Kf_0\oplus \Kf_1$ by:
\begin{eqnarray*}
x\cdot y&=& xy, \quad x,y \in \Kf_0,\\
x\cdot a &=& \textstyle \frac{1}{2} x\star a, \quad x\in \Kf_0,\, a \in \Kf_1,\\
a\cdot x &=&\textstyle \half x\star a, \quad x\in \Kf_0,\, a \in \Kf_1,\\
a\cdot b &=& <a,b>, \quad a,b\in \Kf_1.
\end{eqnarray*}

\subsection{Kaplansky-Jordan superalgebras}
Kaplansky superalgebras were studied by McCrimmon \cite{McC}.
The axioms of Lie antialgebras appear separately in different places of his work. For instance,  the McCrimmon's Theorem 1.2
states:
\textit{
A half-unital superalgebra is a Jordan superalgebra if and only if it satisfies the identities (LA1'), (LA2) and (LA3).
}
Indeed, using the half-unit one can check that the quartic identities (SJ2) of Jordan superalgebras imply the cubic identities (LA1'), (LA2) and (LA3). The proof of the ``if'' part of the theorem is more involved.

From Theorem 1.2 and its proof and Theorem 3.7 of \cite{McC}, we can immediately establish the following relations between the above classes of algebras.

\begin{enumerate}

\item[\textbullet]  Every (weak) Lie antialgebra is a Jordan superalgebra.
\item[\textbullet]  Every Jordan superalgebra which is half-unital and which has an associative even part is a weak Lie antialgebra.
\item[\textbullet]  A Kaplansky superalgebra is a Jordan superalgebra if and only if it is a Lie antialgebra.
\item[\textbullet]  Half-unital Lie antialgebras are exactly Kaplansky-Jordan superalgebras.
\end{enumerate}

\subsection{Two main examples}
The main examples of Kaplansky-Jordan-Ovsienko algebras are:

\begin{enumerate}
\item[(a)]  
\textit{The tiny Kaplansky algebra.} It is a 3-dimensional algebra $K_3$ (denoted by $\asl_2$ in \cite{Ovsienko}) with basis vectors $\e$ (even) and $a,b$ (odd) satisfying:
\begin{equation}
\label{aslA}
\begin{array}{llllllllll}
\e\cdot\e&=&\e,&&&&&&&
\\[4pt]
\e\cdot{}a&=&a\cdot{}\e&=& \half \,a,&
\qquad
\e\cdot{}b&=&b\cdot{}\e&=&\half \,b,
\\[5pt]
a\cdot{}b&=&-b\cdot{}a&=&\half \,\e,&
\qquad
a\cdot{}a&=&b\cdot{}b&=&0.
\end{array}
\end{equation}
\item[(b)]  \textit{The full derivations algebra}, $\A \K(1)$, (called the conformal Lie antialgebra in \cite{Ovsienko}). 
This algebra is generated by even elements $\{\e_n, n\in \Z\}$ and odd elements $\{a_i, \, i\in \Z+\half\}$ satisfying
\begin{equation}
\begin{array}{rcl}
\e_n\cdot \e_m&=& \e_{n+m}\\[5pt]
\e_n\cdot a_i&=& \half a_{n+i}\\[5pt]
a_i\cdot a_j&=& \half (i-j) \e_{i+j}.
\end{array}
\end{equation}

\end{enumerate}

\noindent
The algebra $K_3$ is the unique commutative superalgebra such that $\Der(K_3)=\osp(1|2)$, see \cite{Ovsienko}. Similarly, $\Der(\A \K (1))=\K(1)$, where $\K(1)$ is the conformal Lie superalgebra,  but the uniqueness is still a conjecture.

\subsection{Remark on the axioms of Lie antialgebra}\label{assoc}
Unlike Kaplansky algebras, the even part of Lie antialgebras is not required to contain a unit. 

One can always add a unit $\e$ in $\ga_0$. The algebra $\ga'=(\Kb \e \oplus \ga_0)\oplus \ga_1$ obtained from a non-unital Lie algebra $\ga=\ga_0\oplus \ga_1$ by extending the multiplication to 
$$
\e\cdot a=a, \;a\in \ga_0, \quad \e\cdot a=\textstyle \half a, \;a \in \ga_1,
$$
is still a Lie antialgebra.

The condition of associativity of the even part of Lie antialgebras seems to be natural due to the following proposition.

\begin{prop}\label{axiomassoc}
If the even part $\ga_0$ of a Lie antialgebra $\ga$ is
generated by the odd part $\ga_1$, then the property of associativity 
(LA0) follows from (LA1)--(LA3).
\end{prop}

\begin{proof}
Assume the axioms (LA1)--(LA3).
Using the commutativity, it is equivalent to show
 that $(xy)z=(xz)y$, for all $x, y, z \in a_0$. 

Write $x=ab$ with $a,b\in a_1$.
\begin{eqnarray*}
((ab)y)z&=&((ay)b)z+ (a(by))z\\
&=&((ay)z)b+(ay)(bz)+(az)(by)+a((by)z)\\
&=&\half(a(yz))b+\half a(b(yz))+(ay)(bz)+(az)(by)\\
&=&\half(ab)(yz)+(ay)(bz)+(az)(by)
\end{eqnarray*}
The last expression is symmetric in $y,z$. Therefore, we deduce
$((ab)y)z=((ab)z)y$
and hence the associativity.
\end{proof}

A similar statement holds for Lie superalgebras generated by odd elements. This fact will be established further in the proof of Theorem \ref{thmbrkt}.

\subsection{Representations and modules} 
Any associative superalgebra $\Af$ is a Jordan superalgebra with respect to the product given by so-called \textit{anti-commutator} or \textit{Jordan superproduct}:
\begin{equation}
[a,b]_+=\textstyle \half \big(\;a\cdot b+(-1)^{\bar{a}\bar{b}}b\cdot a\;\big)
\end{equation}
for all homogeneous elements $a,b \in \Af$.

Recall that a Jordan superalgebra $\Jf$ is called \textit{special} if there exist an associative algebra $\Af$ and a faithful homomorphism from $\Jf$ into $(\Af, [\, ,\, ]_+ )$.

A (bi)module over a Jordan superalgebra $\Jf$ is a vector space $V$ together with left and right actions such that the split null extensions $\Jf\oplus V$ and $V\oplus \Jf$ are Jordan superalgebras.

Lie antialgebras are particular cases of Jordan superalgebras. However, we will adopt a slightly different definition of specialization and module in the particular case of Lie antialgebras.  It turns out that a nice theory can be developed after making this change. To avoid any confusion we will be talking of LA-modules and LA-representations. The definitions are the following.

\begin{defn}\cite{Ovsienko}
Let $\ga$ be an arbitrary Lie antialgebra.

(1)
An \textit{LA-module} of $\ga$ is a $\Z_2$-graded vector space $V$ together with  an even linear map
$
\rho:\ga\to\End(V),
$
such that the direct sum $\ga\oplus{}V$ equipped
with the product
\begin{equation}\label{antimod}
(a+v)\cdot(b+w)
=
a\cdot{}b\, +\; \left(
\r(a)\,(w)+(-1)^{\bar{b}\bar{v}}\,\r(b)\,(v) \right),
\end{equation}
where $a,b\in\ga$ and $v,w\in{}V$ are homogeneous elements,
is again a Lie antialgebra.

(2)
An \textit{LA-representation} of $\ga$ is a $\Z_2$-graded vector space $V$ together with  an even linear map
$
\rho:\ga\to\End(V),
$
such that 
\begin{equation}
\label{Firc}
\rho(a\cdot b)=[\rho(a),\rho(b)]_+,
\end{equation}
for all elements $a,b$ in $\ga$, and
\begin{equation}
\label{Secc}
\rho(x)\rho(y)=\rho(y)\rho(x),
\end{equation}
for all even elements $x, y$ in $\ga_0$.
\end{defn}

\begin{rem}
In other words, an LA-representation is a specialization in the usual sense of Jordan superalgebras with
the additional condition of commutativity (\ref{Secc}).
This condition is crucial to link the representations of Lie antialgebras to representations of Lie superalgebras. This condition was first assumed in \cite{Ovsienko}; this is the main difference from the  traditional works on Jordan superalgebras.
\end{rem}

The modules of the tiny Kaplansky algebra $K_3$, considered as a Jordan superalgebra, were studied and classified in \cite{MZ1,MZ2,T1}. The LA-representations of $K_3$ were studied independently in \cite{MG}.

\subsection*{Example: modules over $K_3$}
The classification of finite-dimensional irreducible LA-modules over $K_3$ can be deduced from \cite{MZ2}. 
Let us consider the adjoint module over $K_3$, $i.e$ the 3-dimensional vector space $$V_{\ad}:=\Kb< w> \oplus\;\Kb<u,v> ,$$ 
together with the following action:
\begin{equation}
\begin{array}{lcllcllcl}
\e\cdot v&=&\textstyle \half v\;, \qquad &\e \cdot w&=& w\,, \qquad &\e\cdot u&=&\half u,\\[6pt]
a\cdot v&=& w\;, \quad &a \cdot w&=& u\, ,\quad& a\cdot u&=&0,\\[6pt]
b \cdot v&=&0\;, \quad &b \cdot w&=&\textstyle \frac{1}{4}v\,, \qquad &b\cdot u&=&-\frac{1}{4} w.\\
\end{array}
\end{equation}

\begin{prop} 
Up to isomorphism, $V_\ad$ is the only non-trivial finite dimensional irreducible LA-module of $K_3$.
\end{prop}

\begin{proof} One can easily check that the property \eqref{antimod} in the definition of LA-module holds for  $V_\ad$. From Corollary 2.2 in \cite{MZ2} one knows that the finite dimensional irreducible modules of $K_3$ are the modules $Irr (\sigma,\half, m)$, $\s \in \{0,1\}$, $m\in \Z_{\geq 0}$, in which there exists a vector $v$ of parity $\s$ satisfying:
\begin{equation}\label{modIrr}
\begin{array}{rcl}
\e \cdot v&=&\textstyle \half\, v\\[5pt]
b\cdot v&=& 0\\[5pt]
b \cdot(a \cdot v)&=&\textstyle \frac{m}{4}\,v
\end{array}
\end{equation}
Suppose $Irr (\sigma,\half, m)$ is also an  LA-module. One necessarily has $\sigma=1$. If not the condition $\e \cdot v=\textstyle \half\, v$ leads to a contradiction with the associativity of even elements. Consider the odd element $v \in Irr(1,\half, m)$ satisfying \eqref{modIrr}.
In one hand, using the axiom (LA3) one has
\begin{eqnarray*}
b \cdot (a\cdot v)&=& -a\cdot (v \cdot b)-v\cdot (b\cdot a)\\
&=& \textstyle-v(-\half \e)= \frac{1}{4}\,v
\end{eqnarray*}
But in the other hand, by \eqref{modIrr} one should have 
$$b \cdot(a \cdot v)=\textstyle \frac{m}{4}\,v.$$
This is possible if and only if $m=1$. And in this case $Irr(1,\half, 1)\simeq V_\ad$.
\end{proof}

\section{Universal enveloping algebra.}\label{UEA}

\subsection{Notations} Let $\ga$ be a vector space. We denote by $\T(\ga)$, resp. $\Sc(\ga)$, the universal tensor, resp. symmetric, algebra over $\ga$. Recall briefly, as a vector space 
$$\T(\ga):=\bigoplus_{n\geq 0} \ga ^{\otimes n}.$$
The product of $\T(\ga)$ is given by
$$a_{1}\otimes \cdots \otimes a_{n} \;\cdot\;  b_{1}\otimes \cdots \otimes b_{m} =a_{1}\otimes \cdots \otimes a_{n} \otimes  b_{1}\otimes \cdots \otimes b_{m} $$
on $\ga^{\otimes n}\times \ga^{\otimes m}$ into $\ga^{\otimes n+m}$, and extended by bilinearity on  $\T(\ga)\times  \T(\ga)$.
The algebra $\T(\ga)$ is naturally $\N$-graded. For all $n\in \N$ we denote by 
$$  \T^n(\ga):=   \ga^{\otimes n} $$
the homogeneous component of degree $n$.  

The symmetric algebra is the following quotient
$$\Sc(\ga):=\T(\ga)/<a\otimes b-b\otimes a, \; a,b\in\ga>.$$
We denote by $a\odot b $ the class of $a\otimes b$ in $\Sc(\ga)$. We denote by $\Sc^n(\ga)$ the homogeneous component of degree $n$, that is the image of $\T^n(\ga)$ in $\Sc(\ga)$.

In addition, if $\mathfrak{a}=\ga_0\oplus \ga_1$ is a $\Z_2$-graded vector space, then $\T(\ga)$ is also a $\Z_2$-graded algebra with respect to the following grading
$$\T(\ga)_i:=\bigoplus_{n\geq 0}  \;\bigoplus _{i_1+ \cdots+ i_n=i}\ga _{i_1}\otimes \cdots \otimes \ga_{i_n}, \quad i=0,1.$$
This grading also induces a $\Z_2$-grading on the algebra $\Sc(\ga)$.

\subsection{Definitions} 
Given a Lie antialgebra $\mathfrak{a}=\ga_0\oplus \ga_1$, one associates an associative $\Z_2$-graded algebra $\U(\ga)$ called the \textit{universal enveloping algebra} of $\ga$. The construction of $\U(\ga)$ is as follows. 

Consider the tensor algebra $\T(\ga)$ and denote by $\Pc$ the subset
\begin{equation}
\label{RelaT}
\left\{
\begin{array}{l}
\half \big(a\otimes b -b\otimes a\big) - ab\;, \\[4pt]
\half \big(a\otimes x + x\otimes a\big) - ax\;,\\[4pt]
 x\otimes y - xy,
 \end{array}
 \right.
 \end{equation}
 where $ a,b \in \ga_1$ and $x,y \in \ga_0$.
We define
$$\U(\ga):=\T(\ga)/ <\Pc>,$$
where  $<\Pc>$ denotes the two-sided ideal generated by $\Pc$.


 We introduce the canonical projection
 $$\pi: \T(\ga)\twoheadrightarrow \U(\ga).$$
 When no confusion occurs, we use the same notation for the element $x_1\otimes x_2\otimes \cdots \otimes x_n$ in $\T(\ga)$ and its class in $\U(\ga)$.
 
 We also introduce the canonical embedding of $\ga$ into $\U(\ga)$:
 \begin{equation}\label{defiota}
 \iota:\;\ga \hookrightarrow \ga^{\otimes 1}=\T^1(\ga) \hookrightarrow \T(\ga)\twoheadrightarrow \U(\ga).
 \end{equation}
 
\subsection{Universal property}
The couple $(\U(\ga),\iota)$ satisfies a universal property.

\begin{prop}
Let $\A$ be a $\Z_2$-graded associative algebra and $\phi$ be a linear morphism $\phi: \ga \to \mathcal{A}$ satisfying 
$$
\phi(xy)=\half\big(\phi(x)\phi(y)+(-1)^{\bar{x}\bar{y}}\phi(y)\phi(x)\big),
$$
for all $x,y \in \ga_0 \cup \ga_1$ and $\phi(x)\phi(y)=\phi(y)\phi(x)$ for all  $x,y \in \ga_0$. There exists a unique morphism of algebras $\phi':\U(\ga) \to \mathcal{A}$ such that the following diagram is commutative:

\medskip
 \SelectTips{eu}{12}%
\xymatrix{
   &&&&&& \U(\ga) \ar[dr]^{\phi'}&\\
   &&&&& \ga  \ar[rr]^\phi  \ar[ru]^{\iota}&& \mathcal{A} 
  }
\medskip
\end{prop}

\begin{proof}
The universal property of $\T(\ga)$ gives an homomorphism of algebra $\Phi$ from $\T(\ga)$ to $\A$ such that $\Phi_{|\ga}=\phi$. Since $\mathcal P \subset \Ker(\Phi)$, the morphism  $\Phi$ induces a morphism $\phi'$ from $\U(\ga)$ to $\A$ such that $\phi=\phi' \circ \iota$.
\medskip

 \SelectTips{eu}{12}%
\xymatrix{
 &&&&&\T(\ga) \ar@{->>}[r] \ar@{-->}[rrd]_{\Phi}& \U(\ga) \ar[dr]^{\phi'}&\\
   &&&&& \ga \ar@{^{(}->}[u]   \ar[rr]_\phi  \ar[ru]^{\iota}&& \mathcal{A} 
   }
 \medskip
 \noindent
Since $\U(\ga)$ is generated by the elements $\iota(x)$, $x\in \ga$, the condition $\phi' \circ \iota =\phi$ uniquely determines $\phi'$.
\end{proof}

\subsection{The Poincar\'e-Birkhoff-Witt property}
The natural $\N$-filtration of  $\U(\ga)$  inherited from the natural filtration  of $\T(\ga)$, is given by
$$
 \U_{n}:=  \pi 
 \Big( \bigoplus_{0\leq k \leq n} \ga^{\otimes k}\,\Big), \; n\in \N.
$$
The associated graded algebra is
$$\Gr \,\U:=\bigoplus_{n\geq 0} \U_{n}/ \U_{n-1},$$
with the usual conventions $\U_0=\mathbb{K}$ and $\U_{-1}=\{0\}$.

Consider the algebra
$$G(\ga)= \T(\ga)/<\Rc>,$$
where  $\Rc$ is the set of homogeneous quadratic relations obtained from the relations (\ref{RelaT}) by projection onto $\ga\otimes\ga$.
More precisely,
$$
\Rc=\left\{a\otimes b-b\otimes a\, ,
\quad a\otimes x + x\otimes a\,,
\quad x\otimes y\right\},$$
where $x,y\in \ga_0, a,b\in\ga_1$.

The algebra $\U(\ga)$ satisfies the Poincar\'e-Birkhoff-Witt (PBW) property if 
$$
\Gr \,\U\cong G(\ga).
$$
We will show that, in general, the universal enveloping algebra of
a Lie antialgebra does not satisfy the PBW property.

\begin{prop}
\label{MalPr}
One has an isomorphism of algebra
$$ G(\ga)\simeq\; (\Kb \oplus \ga_0)\otimes \Sc(\ga_1),$$
where the product on  $ (\Kb \oplus \ga_0)\otimes \Sc(\ga_1) $ is given by 
$$
(\lambda +x) \otimes a \;\cdot\; (\mu +y) \otimes b = (\lambda \mu +\lambda x+\mu y)\otimes a\odot b, 
$$
for all $\lambda, \mu \in \Kb$, $x,y \in \ga_0$, $a,b \in \Sc(\ga_1)$.
\end{prop}
\begin{proof}
If $\{x_i, i\in I\}$ is a basis of $\ga_0$ and $\{a_j, j\in J\}$ is a basis of $\ga_1$, where the index set and $J$ is totally ordered, then a basis of $G(\ga)$ is given by the set of monomials,
$$
x_i a_{j_1}\cdots  a_{j_p}, \quad \text{ and }\quad   a_{j_1} \cdots a_{j_p}, 
$$
where $i \in I$, $p\in \N$, and $j_1\leq \cdots \leq j_p$ is an increasing sequence of indices in $J$. The multiplication of basis elements in $G(\ga)$ is
\begin{equation}\label{multGa}
\left\{
\begin{array}{lcl}
x_i a_{j_1} \cdots a_{j_p}\cdot x_{i'} a_{j'_1} \cdots a_{j'_q}&=&0,\\ a_{j_1} \cdots a_{j_p}\cdot x_{i'} a_{j'_1} \cdots a_{j'_q}&=&(-1)^px_{i'} a_{j_1} \cdots a_{j_p} a_{j'_1} \cdots a_{j'_q}.
\end{array}
\right.
\end{equation}
Hence the result.
\end{proof}

\begin{prop}
For a Lie antialgebra $\ga$, the algebra $\Gr \,\U(\ga)$ is not necessarily isomorphic to $G(\ga)$. 
\end{prop}
\begin{proof}
This can be deduced by using a necessary condition given by Braverman and Gaitsgory in \cite{BG}. Following \cite{BG}, denote by $\alpha : \Span _\Kb \Qc \rightarrow \ga$ the linear map
\begin{eqnarray*}
\alpha(a\otimes b-b\otimes a)&=&2\,ab\\
\alpha(a\otimes x + x\otimes a)&=&2\,ax\\
\alpha(x\otimes y)&=&xy,
\end{eqnarray*}
for $a,b\in\ga_1$ and $x,y\in\ga_0$.

The set of relations $\Pc$ then can be described as
$$
\Pc=\{\,q-\alpha (q)\,,\; q\in \Qc\,\}.
$$

\noindent
Let us refer to \textit{BG-conditions} the following two conditions:

\begin{enumerate}
\item[(i)] $\im (\alpha \otimes \id-\id \otimes \alpha )\subset  \Span _\Kb \Qc $ (where $\alpha \otimes \id-\id \otimes \alpha$ is defined on the space $\Span _\Kb (\Qc \otimes \ga) \cap \Span _\Kb (\ga \otimes \Qc)$),
\item[(ii)] $\a \circ  (\alpha \otimes \id-\id \otimes \alpha )=0$.
\end{enumerate}

According to \cite{BG}, if $\U(\ga)$ has the PBW-property, then the BG-conditions are satisfied.

Let us show that in our situation the BG-conditions fail. Indeed, the space 
 $ \Span _\Kb (\Qc \otimes \ga) \cap \Span _\Kb (\ga \otimes \Qc)$ is generated by
 the following elements:
 \begin{equation}\label{span}
 \left\{
\begin{array}{lcl} 
 u_0&=& (x\otimes y)\otimes z\\
 &=&x\otimes (y\otimes z),\\[6pt]
 u_1&=&(a\otimes x+x\otimes a)\otimes y +(x\otimes y)\otimes a\\
 &=& a\otimes (x\otimes y)+x\otimes (a\otimes y+y\otimes a),\\[6pt]
 u_2&=& x\otimes (a\otimes b-b\otimes a)-b\otimes(x\otimes a+a\otimes x)+a\otimes (b\otimes x+x\otimes b)\\
&=&(x\otimes a+a\otimes x)\otimes b-(x\otimes b+b\otimes x)\otimes a+(a\otimes b-b\otimes a)\otimes x,\\[6pt]
u_3&=& (a\otimes b-b\otimes a)\otimes c+(b\otimes c-c\otimes b)\otimes a+(c\otimes a-a\otimes c)\otimes b\\
&=&a\otimes (b\otimes c-c\otimes b)+b\otimes (c\otimes a-a\otimes c)+c \otimes (a\otimes b-b\otimes a),
\end{array}
\right.
\end{equation}
where $a,b,c$ are elements  in $\ga_1$ and $x,y,z$ in $\ga_0$.

One can check that the elements $(\alpha \otimes \id-\id \otimes \alpha )(u_i)$ are not elements of $\Span \Qc$ in general, see Section \ref{KAK1} below for a counterexample.
\end{proof}

\begin{rem}
In the case of $\ga= K_3$, the BG-conditions hold. Indeed, the set of generators described above in \eqref{span} simplifies to
 \begin{equation*}
 \left\{
\begin{array}{lcl} 
 u_0&=& \e \otimes \e \otimes \e\\[6pt]
 u_1&=&(t\otimes \e+\e\otimes t)\otimes \e +(\e \otimes \e)\otimes t,\\
 &=& t\otimes (\e\otimes \e)+\e\otimes (t\otimes \e+\e\otimes t), \quad t=a,b,\\[6pt]
 u_2&=& \e\otimes (a\otimes b-b\otimes a)-b\otimes(\e\otimes a+a\otimes \e)+a\otimes (b\otimes \e+\e\otimes b)\\
 &=&(\e\otimes a+a\otimes \e)\otimes b-(\e\otimes b+b\otimes \e)\otimes a+(a\otimes b-b\otimes a)\otimes \e,\\[6pt]
\end{array}
\right.
\end{equation*}
where $\e, a, b$ are the elements of the standard basis of $K_3$. 

Using the relations in $K_3$, one easily gets,
\begin{eqnarray*}
(\alpha \otimes \id-\id \otimes \alpha )(u_i)&=&0, 
\end{eqnarray*}
for every generators  $u_i$, $i=0,1,2$.
This implies the BG-conditions. 

In \cite{BG}, it is shown that in the case where $G(\ga)$ is of Koszul type the BG-conditions are sufficient to imply the PBW property. We do not know if $G(K_3)$ is of Koszul type so we will show by hand that $\U(K_3)$ satisfies the PBW theorem (see Proposition \ref{PBWK3}).
\end{rem}

\subsection{Example: the enveloping algebra $\U(K_3) $}
The algebra $\U(K_3) $ is the associative algebra generated by two odd elements $A, B$ and one even element $\E$ subject to the relations 
\begin{equation} \label{relUa}
\left \lbrace
\begin{array}{rcl}
AB -BA&=&\E\\ [5pt] 
A\E + \E A&=&A\\ [5pt]
B\E + \E B&=&B\\ [5pt]
 \E^2&=&\E.
\end{array}
\right.
\end{equation}

\noindent
One can easily establish the following useful formulae by induction:
\begin{lem}\label{formulaire}
For any integers $p,q,k,l$, one has
\begin{enumerate}
\item $A^{2p}\E=\E A^{2p}$ and $A^{2p+1}\E=A^{2p+1}-\E A^{2p+1}$,
\item $B^{2q}\E=\E B^{2q}$ and $B^{2q+1}\E=B^{2q+1}-\E B^{2q+1}$,
\item if $k+l$ is even then $A^kB^l\E=\E A^kB^l$,
\item if $k+l$ is odd then $ A^kB^l \E =A^kB^l-\E A^kB^l$, in particular  $ \E A^kB^l \E =0$,
\item $ BA^{2p}=A^{2p}B-p A^{2p-1}$ and  $BA^{2p+1}=A^{2p+1}B-\E A^{2p}-pA^{2p}$,
\item $B^{2q}A=AB^{2q}-qB^{2q-1}$ and $B^{2q+1}A=AB^{2q+1}-\E B^{2q}-qB^{2q}$.
\end{enumerate}
\end{lem}

\begin{prop} The following set of monomials 
$$
\{ \E A^k B^l\, , \; k,l\in \N\} \cup \{  A^k B^l\, , \; k,l\in \N\}
$$
is a basis of $\U(K_3)$.
\end{prop}

\begin{proof}
Using Lemma \ref{formulaire}, one can see that the above set of monomials is spanning the vector space $\U(K_3)$. To see that they are linearly independent we use weight and order functions, and a representation of $\U(K_3)$ in the space $\mathscr{D}$ of differential operators over the algebra $\C[x, \xi]/<\xi^2>$.

Let us introduce the following weight function
$$\wt(A)=1, \quad \wt(\E)=0, \quad \wt(B)=-1.$$
The relations \eqref{relUa} are homogeneous with respect to the weight function. Therefore, linear dependence of the monomials could only occur between monomials of same weight. Let us fix a weight $w$ in $\Z$ and show that the monomials
\begin{equation}
\label{TheSet}
\{ \E A^k B^l\, , \; k-l=w\} \cup \{  A^k B^l\, , \; k-l=w\}
\end{equation}
are linearly independent.

Let us consider the following operator $\Dc$ of order $\half$ in $\mathscr{D}$
$$
\Dc:=\partial_\xi + \xi \partial_x.
$$
It is easy to check (cf. \cite{MG}) that the map defined by
$$
A\mapsto x\Dc,
\qquad
B\mapsto \Dc,
\qquad
\E\mapsto \xi \Dc
$$
is a representation of $\U(K_3)$ into the space $\mathscr{D}$. In this representations, one has
\begin{equation*}
\left\{
\begin{array}{lcl}
\E A^kB^l&\mapsto& \xi x^k\Dc^{k+l+1}+\text{ operators of order }< (k+l+1)/2\\[6pt]
A^kB^l&\mapsto& x^k\Dc^{k+l}+\text{ operators of order }< (k+l)/2.
\end{array}
\right.
\end{equation*}
Therefore, all of the differential operators corresponding to the monomials of the set
(\ref{TheSet})
have different orders.
It follows that these monomials are linearly independent and thus form a basis of $\U(K_3)$.
\end{proof}

\begin{prop}\label{PBWK3}
The enveloping algebra $U(K_3)$ satisfies the PBW property.
\end{prop}

\begin{proof}
Using Lemma \ref{formulaire}, one can see that 
\begin{equation*}
\left\{
\begin{array}{lcl}
\! \E A^kB^l\cdot \E A^{k'}B^{l'}\!&=\!&\!(-1)^{k+l}\E A^{k+k'}B^{l+l'}\!+\text{terms of length}\leq (k+l+k'+l'+1)\\[6pt]
A^kB^l\cdot \E A^{k'}B^{l'}\!&=\!&\! (-1)^{k+l}\E A^{k+k'}B^{l+l'}\!+\text{terms of length}< (k+l+k'+l').
\end{array}
\right.
\end{equation*}
So, in $\Gr\; \U(K_3)$ one gets
\begin{equation*}
\left\{
\begin{array}{lcl}
\E A^kB^l\cdot \E A^{k'}B^{l'}&\equiv&0\\[6pt]
A^kB^l\cdot \E A^{k'}B^{l'}&\equiv& (-1)^{k+l}\E A^{k+k'}B^{l+l'}\end{array}
\right.
\end{equation*}
what is the same multiplication rules as in $G(K_3)$, see \eqref{multGa}.
\end{proof}

\begin{rem}
 The enveloping algebra $\U(\ga)$ endowed with the  Jordan superproduct
$$
[x, y]_+ =\half (x\otimes y + (-1)^{\bar{x}\bar{y}}y\otimes x ).
$$
 is not a Lie antialgebra.
One can notice that in our example the axiom (LA1) fails. Indeed, using Lemma \ref{formulaire} one checks
\begin{eqnarray*}
[\E AB,[ \E B^2, \E B]_+]_+&=& \half [\E AB, \E B^2 \E B+ \E B\E B^2]_+\\
&=& \half [\E AB,  \E B^2 \E B]_+ \\
&=&\quart (\E AB \E B^2 \E B+ \E B^2 \E B \E AB)\\
&=&\quart  \E AB^4,
\end{eqnarray*}
whereas
\begin{eqnarray*}
[[\E AB, \E B^2]_+, \E B]_+&=& \half [\E AB  \E B^2+  \E B^2\E AB\,,\;  \E B]_+\\
&=&\quart (\E AB  \E B^2 \E B+  \E B^2\E AB  \E B + 
\E B \E AB  \E B^2+\E B \E AB)\\
&=&\quart ( \E AB^4 + \E B^2 A B^2)\\
&=&\quart ( \E AB^4  + \E(AB^2-B)B^2)\\
&=& \half \E AB^4 - \quart \E B^3.
\end{eqnarray*}
Thus one gets
$[\E AB,[\E  B, \E B^2]_+]_+\,\not=\frac{1}{2}[[\E AB, \E B^2]_+, \E B]_+$.
\end{rem}

\section{Links between Lie superalgebras and Lie antialgebras}\label{liens}

\subsection{Adjoint Lie superalgebra}\label{adjSuper}

In this section, we study the construction of the Lie superalgebra associated to a Lie antialgebra \cite{Ovsienko}. We provide the missing proofs of \cite{Ovsienko} in Section \ref{ProofSect}.

Given a Lie antialgebra $\ga$, the \textit{adjoint  Lie superalgebra} $\gog_\ga$  is defined as follows.
As a vector space $\gog_\ga=(\gog_\ga)_0\oplus (\gog_\ga)_1$, where
$$
(\gog_\ga)_1:=\ga_1, \quad (\gog_\ga)_0:=\ga_1 \otimes \ga_1/ \mathcal{S}
$$
and $\mathcal{S} $ is the ideal generated by
$$\{a\otimes b-b\otimes a,\; ax\otimes b-a\otimes bx\,| \; a,b\in \ga_1, x \in \ga_0\}.$$
We denote by $a\odot b$ the image of $a\otimes b$ in $(\gog_\ga)_0$. Therefore, we have the following useful relations in $(\gog_\ga)_0$:
\begin{equation}\label{relequiv}
\left\{
\begin{array}{rcl}
a\odot b& = &b\odot a,\\
ax\odot b& = &a\odot bx \;= \;b\odot ax\; =\;bx\odot a, \quad \;a,b\in \ga_1, x \in \ga_0 .\\
\end{array}
\right.
\end{equation}

\noindent
The Lie bracket on $\gog_\ga$ is given by:
\begin{equation}\label{superbracket}
\begin{array}{rcl}
[a,b]&=& a\odot b,\\[5pt]
[a\odot b, c]=-[c,a\odot b]&=&a(bc)+b(ac),\\[5pt]
[a\odot b, c\odot d]&=& 2\,a(bc)\odot d + 2\, b(ad)\odot c,
\end{array}
\end{equation}
where $a,b,c$ and $d$ are elements of $(\gog_\ga)_1=\ga_1$.

\begin{ex}
In the $\ga=K_3$ case the above construction leads to $\gog_\ga=\mathrm{osp}(1|2)$. 
More precisely,
$$
(\gog_\ga)_1=<a,b>, \quad (\gog_\ga)_0=<a\odot a, \,a\odot b, \,b\odot b>.
$$
Rescale the basis elements as follows:
$$E:= 2 a\odot a,\; F:=-2 b\odot b, \; H:=- 4 a\odot b,\; A:=2a,\; B:=2b.$$
Then using  the definitions of the bracket and the relations in $K_3$ one checks
\begin{equation}
\label{commutations}
\begin{array}{lll}
[H, E]=2E,  &[H,F]=-2F, &[E,F]=H,\\[5pt]
[H, A]=A,  &[E,A]=0,  &[F,A]=B,\\[5pt]
[H, B]=-B,            &[E,B]=A,              &[F,B]=0,\\[5pt]
[A, B]=-H, &[A, A]= 2E, & [B, B]= -2F.
\end{array}
\end{equation}
The above system is a presentation of $\osp(1|2)$.

\end{ex}

\begin{rem}
This construction differs from the Kantor-Koecher-Tits construction. For instance, the Lie superalgebra associated to $K_3$ through KKT is $\mathrm{psl}(2,2)$, see~\cite{Kac}. In general, if one applies KKT to a Jordan superalgebra $\ga=\ga_0\oplus \ga_1$,  the odd subspace of the resulting Lie superalgebra is much bigger than $\ga_1$.
\end{rem}

In \cite{Ovsienko}, the following statements are formulated.

\begin{prop}\label{propbrkt}
The definition of the bracket on $\gog_\ga$ given in \eqref{superbracket} is compatible with the relations \eqref{relequiv} in $(\gog_\ga)_0$.
\end{prop}

\begin{thm}\label{thmbrkt}
The bracket given in \eqref{superbracket} is a Lie superbracket on $\gog_\ga$. 
\end{thm}

\noindent
We give the direct proofs of these statements in Section \ref{ProofSect}.

 \subsection{Relations between the enveloping algebras}

The following proposition relates the universal enveloping algebra of an arbitrary Lie antialgebra $\ga$ to the  universal enveloping algebra of the adjoint Lie superalgebra  $\gog_\ga$.
\begin{thm}
\label{quotient}
Let $\ga $ be a Lie antialgebra such that the odd part $\ga_1$ spans the even part $\ga_0$. If $\gog:=\gog_\ga$ is the Lie superalgebra associated to $\ga$ then,
$$\U(\ga)\simeq \U(\gog)/\mathcal{I}_\ga$$
for some ideal $\mathcal{I}_\ga$ of $ \U(\gog)$.
 \end{thm}
 
 \begin{proof}
 Recall that $U(\gog)$ is the quotient of the universal associative algebra
  $$T(\gog):=\bigoplus_{n\geq 0} \gog^{\otimes n}$$ 
 by the the 2-sided ideal 
 $$\mathcal{J}: =<\half (x\otimes y -(-1)^{\bar x\bar y}y\otimes x) - [x,y]>\;, \quad x,y \text{ homogeneous in } \gog. $$
 The algebra $T(\gog)$ is spanned by the elements of $\gog$. We define an homomorphism of algebra 
 $$\pi: T(\gog) \rightarrow \U(\ga)$$
 by setting
  \begin{equation}\label{defpi}
\pi(a)=a, \quad \pi(a\odot b)=\half (a\otimes b + b\otimes a)
\end{equation}
 for all odd elements $a$ and $b$ (recall they can be viewed as elements of $\ga_1$ or $\gog_1$).
 The application $\pi$ is well defined. Indeed, one can check that $ \pi(ax\odot b)= \pi(a\odot bx)$ for all $a, b \in \ga_1$ and $x\in \ga_0$:
 
 \begin{equation}\label{pidef}
 \begin{array}{lcl}
\pi(ax\odot b)- \pi(a\odot bx)
&=& \half(ax\otimes b + b\otimes ax-a\otimes bx - bx\otimes a)\\[4pt]
&=& \half (\half (a\otimes x+x\otimes a)\otimes b + b\otimes \half (a\otimes x+x\otimes a)\\[4pt]
&&-a\otimes \half (b\otimes x+x\otimes b) - \half (b\otimes x+x\otimes b)\otimes a)\\[4pt]
&=&\quart (x\otimes a\otimes b+ b\otimes a\otimes x-a\otimes b\otimes x-x\otimes b\otimes a)\\[4pt]
&=&\quart (x\otimes (a\otimes b-b\otimes a)-(a\otimes b-b\otimes a)\otimes x)\\[4pt]
&=&\half (x\otimes ab-ab\otimes x)\\[4pt]
&=&0.
 \end{array}
 \end{equation}
 
 \begin{lem}
 The morphism $\pi$ is surjective.
 \end{lem}

\begin{proof}
We are under the assumption that the odd part of $\ga$ spans the algebra. By definition of $\pi$ all the odd elements are reached.
\end{proof}

 \begin{lem}\label{lemquotient}
 With the above notations, one has
 $$\mathcal{J}\subset \Ker ( \pi).$$
 \end{lem}
 \begin{proof}
 One needs to check 
 \begin{equation}\label{kerpigeneral}
 \pi\big(\half(x\otimes y -(-1)^{\bar x \bar y}y\otimes x )\big)=\pi([x,y])
 \end{equation}
 for all homogeneous elements of $\gog$.
 The more involved case is the one where $x$ and $y$ are both even elements \textit{i.e}
 \begin{equation}\label{kerpi}
 \pi(a\odot b \otimes c\odot d -c\odot d \otimes a \odot b)=
2 \;\pi([a\odot b, c \odot d]).
 \end{equation}
 Let us expand the the left hand side (LHS) of \eqref{kerpi} using the definition of $\pi$, we get :
 \begin{eqnarray*}
 \mathrm{LHS}&=&\textstyle \half(a\otimes b+ b \otimes a)\otimes \half(c \otimes d +d\otimes c)-\half(c \otimes d +d\otimes c)\otimes \half(a\otimes b+ b \otimes a)\\
 &=& \quart(a\otimes b \otimes c \otimes d+ a\otimes b \otimes d \otimes c+  b\otimes a \otimes c \otimes d+ b\otimes a \otimes d \otimes c\\
 &&  \quad -c\otimes d \otimes a \otimes b -c\otimes d \otimes b \otimes a -d\otimes c \otimes a \otimes b-d\otimes c \otimes b \otimes a )  .
 \end{eqnarray*}
Now let us expand the right hand side (RHS)  of \eqref{kerpi} using the equivalent expression \eqref{supbkt2} of the bracket and then the definition of $\pi$. We get:
 \begin{eqnarray*}
  \mathrm{RHS}&=& 2\pi(a(bc)\odot d -c(da)\odot b+b(ac)\odot d-c(db)\odot a)\\
&=&  a(bc)\otimes d +d \otimes a(bc)-c(da)\otimes b- b\otimes c(da)+b(ac)\otimes d+d\otimes b(ac)\\
&&\quad-c(db)\otimes a-a\otimes c(db).
  \end{eqnarray*}
Now using the relations in $\U(\ga)$ we can rewrite each terms as tensor products of four elements. For instance 
 \begin{eqnarray*}
a(bc)\otimes d &=&\half (a\otimes bc+bc\otimes a)\otimes d\\
&=&\quart (a\otimes b \otimes c -a\otimes c \otimes b)\otimes d.
  \end{eqnarray*}
Hence we write
\begin{eqnarray*}
\mathrm{4.RHS}&=& (a\otimes b \otimes c -a\otimes c \otimes b + b\otimes c \otimes a -c\otimes b \otimes a)\otimes d \\
&&\quad + d \otimes (a\otimes b \otimes c -a\otimes c \otimes b + b\otimes c \otimes a -c\otimes b \otimes a)\\
&&\quad -(c\otimes d \otimes a -c\otimes a \otimes d + d\otimes a \otimes c -a\otimes d \otimes c)\otimes b \\
&&\quad -b \otimes (c\otimes d \otimes a -c\otimes a \otimes d + d\otimes a \otimes c -a\otimes d \otimes c)\\
&&\quad +(b\otimes a \otimes c -b\otimes c \otimes a + a\otimes c \otimes b -c\otimes a \otimes b)\otimes d \\
&&\quad + d \otimes (b\otimes a \otimes c -b\otimes c \otimes a + a\otimes c \otimes b -c\otimes a \otimes b)\\
&&\quad -(c\otimes d \otimes b -c\otimes b \otimes d + d\otimes b \otimes c -b\otimes d \otimes c)\otimes a\\
&&\quad - a \otimes (c\otimes d \otimes b -c\otimes b \otimes d + d\otimes b \otimes c -b\otimes d \otimes c).
\end{eqnarray*}
Some of the above terms obviously cancel. The remaining terms can be reorganized as follows:
\begin{eqnarray*}
\mathrm{4.RHS}&=& \mathrm{4.LHS}+(b \otimes c -c\otimes b)\otimes (a\otimes d -d\otimes a)- (a\otimes d -d\otimes a) \otimes (b \otimes c -c\otimes b)\\
&& \quad -(b \otimes d -d\otimes b)\otimes (a\otimes c -c\otimes a)- (a\otimes c -c\otimes a) \otimes (b \otimes d -d\otimes b).
\end{eqnarray*}
Using relations in $\U(\ga)$ this simplify to 
\begin{eqnarray*}
\mathrm{RHS}&=& \mathrm{LHS}+bc\otimes ad-ad\otimes bc -ac \otimes bd +bd \otimes ac.
\end{eqnarray*}
But, in $\U(\ga) $ even elements of order 1 commute so the last terms cancel and we have established \eqref{kerpigeneral} in the case where $x$ and $y$ are two even elements. The other cases are straightforward.
\end{proof}

Since $\mathcal{J}\subset \Ker ( \pi)$, the homomorphism $\pi$ induces a surjective homomorphism $\tilde{\pi}$ from $\U(\gog)$ to $\U(\ga)$. One has $\U(\ga)\simeq \U(\gog)/\Ker (\tilde{\pi})$.

Theorem \ref{quotient} is proved.
 \end{proof}
 
\begin{ex}  In the case $\ga=K_3$ it was already noticed in \cite{MG} that 
$$\U(K_3)=\U(\osp(1|2))\:/ <C>,$$
where $C$ is the usual Casimir element of $\U(\osp(1|2))$. Using the generators $E,F$ and $H$ and $A, B$ of  $\U(\osp(1|2))$ subject to the relations \eqref{commutations} for the Lie superbracket, we express $C$ as
\begin{eqnarray*}
C&=&\textstyle E{}F+F{}E+
\frac{1}{2}
H^2+ \half \left(A{}B-B{}A
\right).
\end{eqnarray*}

An alternative presentation of the ideal is also given in terms of ghost Casimir (see \cite{MG} for more details). The ghost casimir $\cc$ of $\U(\osp(1|2))$ can be expressed as 
$$\cc=AB-BA-\textstyle \half.$$
It is the unique element invariant for the twisted adjoint action
\begin{equation}
\label{Twist}
\widetilde{\ad}_XY:=
XY-(-1)^{\bar X( \bar Y+1)}\,YX.
\end{equation}
for $X, Y$ in $\U(\osp(1|2))$.

Equivalently, one has
$$\U(K_3)=\U(\osp(1|2))\:/ <\textstyle \cc^2-\frac{1}{4}>.$$
This presents $\U(K_3)$ as a generalized Weyl algebra as given in \cite{Bavula}.
\end{ex}

 \subsection{Extension of representations}
In \cite{Ovsienko}, Ovsienko claimed that any representation of a Lie antialgebra $\ga$ can be extended to a representation of the Lie superalgebra $\gog_\ga$. This property can be viewed as a consequence of Theorem \ref{quotient}.

\begin{cor}\label{proprep} Let $(V, \rho)$ be a representation of a Lie antialgebra $\ga$. We define a map  $\tilde{\rho}$ on $\gog_a$ by setting for all $a,b\in\ga_1$:
\begin{equation}
\begin{array}{lcll}
\tilde{\rho}(a)&:=&\rho(a) \\
\tilde{\rho}(a\odot b)&:=&\half(\,\rho(a)\rho(b)+\rho(b)\rho(a)\,).
\end{array}
\end{equation}
The map $\tilde{\rho}$ 
is a representation of $\gog_\ga$.
\end{cor}

\begin{proof}

By universal property of $\U(\ga)$, the representation $\rho$ induces an algebra homomorphism $\rho'$ from $\U(\ga)$ to $\End(V)$. The map $\tilde \rho ':= \rho'\circ \pi$, where $\pi$ is the surjection defined in  \eqref{pidef}, is a homomorphism from $\U(\gog_\ga)$ to $\End(V)$. The restriction of $\tilde \rho '$ to $\gog_\ga$ is a representation of $\gog_\ga$, namely  $\tilde \rho$.
 
\medskip

 \SelectTips{eu}{12}%
\xymatrix{
 &&& \U(\ga) \ar[dr]^{\rho'}&&\U(\gog_\ga)\ar@{->>}[ll]_\pi\ar[dl]_{\tilde \rho '}\ar@{<-^{)}}[dr]\\
   && \ga   \ar[rr]_\rho \ar@{^{(}->}[ru]^{\iota}&& \End(V)&&\gog_\ga\ar@{-->}[ll]^{\tilde \rho}
   }
 \medskip 

\noindent
Hence the result.
\end{proof}

\section{The conformal Lie antialgebra $\A\K(1)$} \label{KAK1}

An interesting example of Lie antialgebra which plays an important role in \cite{Ovsienko} is the infinite-dimensional algebra $\A \K(1)$ called the conformal Lie antialgebra.
This algebra is generated by even elements $\{\e_n, n\in \Z\}$ and odd elements $\{a_i, \, i\in \Z+\half\}$ satisfying
\begin{equation}
\left\{
\begin{array}{ccl}
\e_n\, \e_m&=& \e_{n+m}\\[5pt]
\e_n\, a_i&=& \half a_{n+i}\\[5pt]
a_i\, a_j&=& \half (j-i) \e_{i+j}.
\end{array}
\right.
\end{equation}
It contains infinitely many subalgebras isomorphic to $K_3$.
Note that (a slightly different version of) $\A \K(1)$ was considered in \cite{McC} under the name of full derivation algebra.

\subsection{The universal enveloping algebra $\U(\A\K(1))$} 
The universal enveloping algebra $\U(\A\K(1))$ is the associative algebra generated by  $\{\E_n, n\in \Z,\; A_i, i\in \Z+\half\}$ and the relations:
\begin{equation}\label{relAK1}
\left\{
\begin{array}{rcl}
\E_n\E_m&=&\E_{n+m}\\[5pt]
\E_nA_i+A_i\E_n&=& A_{n+i}\\[5pt]
A_iA_j-A_jA_i&=&(j-i)\E_{i+j}\;.
\end{array}
\right.
\end{equation}

A remarkable additional relation is satisfied in $\U(\A\K(1))$. 

\begin{prop}\label{remcontrex}
One has
\begin{equation}
\label{remaddrel}
A_iA_j+A_jA_i=A_kA_l+A_lA_k,
\end{equation}
for all $i+j=k+l$.
\end{prop}
\begin{proof}
Fix $i,j,k,l$ such that $i+j=k+l$, one gets
\begin{eqnarray*}
A_iA_j+A_jA_i &=&(\E_{i-k}A_k+A_k\E_{i-k})A_j+A_j \, (\E_{i-k}A_k+A_k\E_{i-k})\\
&=& \E_{i-k}A_kA_j+A_k\E_{i-k}A_j+A_j \E_{i-k}A_k+A_jA_k\E_{i-k}\\
&=& \E_{i-k}A_kA_j+A_k(\E_{i-k}A_j+A_j\E_{i-k})-A_kA_j \E_{i-k}\\
&&\quad +(A_j \E_{i-k}+\E_{i-k}A_j)A_k-\E_{i-k}A_jA_k+A_jA_k\E_{i-k}\\
&=& \E_{i-k}(A_kA_j-A_jA_k) + A_kA_l+ A_lA_k-(A_kA_j-A_jA_k)\E_{i-k}\\
&=&(k-j) \E_{i-k}\E_{k+j}+ A_kA_l+ A_lA_k-(k-j)\E _{k+j} \E_{i-k}\\
&=& A_kA_l+A_lA_k.
\end{eqnarray*}
Hence the result.
\end{proof}

\begin{cor}
The algebra $\U(\A\K(1))$ does not satisfy the PBW property. 
\end{cor}
\begin{proof}
The identity (\ref{remaddrel}) can be written as
$$
2A_iA_j-(i-j)\E_{i+j}=2A_kA_l-(k-l)\E_{k+l}, \quad \forall \; i+j=k+l.
$$
Therefore, in $\Gr(\U(\A\K(1))$ one has $A_iA_j=A_kA_l$ for all $i+j=k+l$ and this is not true in $G(\A\K(1))=(\Kb \oplus <\E_n, n\in \Z>)\otimes \Sc(<A_i, i\in \half+\Z>)$.
\end{proof}

\begin{rem}
The canonical mapping $\iota:\A\K(1)\rightarrow \U(\A\K(1))$ introduced in \eqref{defiota} is injective. 
Indeed, if one associates the following weights to the generators
$$\wt(A_i)=i, \quad \wt(\E_n)=n,$$
one can see that the relations in $\U(\A\K(1))$ are homogeneous with respect to this weight function.
Therefore, the elements $\{A_i, i\in \half +\Z, \E_n, n\in \Z\}$ are linearly independent.
\end{rem}


\subsection{Adjoint Lie superalgebra}
One can check (cf. \cite{Ovsienko}) that the adjoint Lie super algebra of $\A \K(1)$ is the conformal (or centerless Neveu-Schwartz) Lie algebra $\K(1)$ generated by
$$
\textstyle
\left\{
x_n,\;n\in\Z;
\qquad
a_i,\;i\in\Z+\half
\right\}
$$
with the following commutation relations
\begin{equation}
\label{CAlgRel}
\left\{
\begin{array}{lcl}
\left[
x_n,x_m
\right] &=&
\half \left(m-n\right)x_{n+m},\\[8pt]
\left[
x_n,a_j
\right] &=&\half
\left(i-\frac{n}{2}\right)a_{n+i},\\[8pt]
\left[
a_i,a_j
\right] &=&
\,x_{i+j}.
\end{array}
\right.
\end{equation}

The corresponding universal enveloping algebra 
$\U(\K(1))$ is the associative algebra generated by  $\{X_n, n\in \Z,\; A_i, i\in \Z+\half\}$ and the relations:
\begin{equation}\label{relK1}
\left\{
\begin{array}{rcl}
X_nX_m-X_mX_n&=& (m-n)X_{n+m}\\[7pt]
X_nA_i-A_iX_n&=&(i-\frac{n}{2}) A_{n+i}\\[7pt]
A_iA_j+A_jA_i&=& 2X_{i+j}.
\end{array}
\right.
\end{equation}

\begin{prop}
One has
\begin{equation}\label{isoAK1}
U(\A\K(1))\cong
\U(\K(1))/
\left<
\begin{array}{ll}
&A_{n-i}A_i-A_iA_{n-i}=(n-2i)\E_n\, , \\[6pt]
& \E_n\E_m=\E_{n+m}, \;n,m\in \Z, i\in \half +\Z
\end{array}
\right>.
\end{equation}
\end{prop}
\begin{proof}
By Lemma \ref{lemquotient}, $\pi(A_i)=A_i$ defines a homomorphism from $\U(\K(1))$ to $\U(\A\K(1))$. In other words, the relations \eqref{relK1} are satisfied in $\U(\A\K(1))$ (see remark \ref{remcontrex} and Example \ref{exrelation} below for direct verifications). 
 
Conversely, in $\U(\K(1))$, let us denote for all $n\in \Z$,
$$
\E_n:=\frac{1}{n-1}( \textstyle A_{n-\half}A_{\half}-A_{\half} A_{n-\half}).
$$
The result follows.
\end{proof}

\begin{ex}\label{exrelation}
The isomorphism \eqref{isoAK1} can be checked by direct computations. For instance, one can use the following relations
\begin{eqnarray*}
X_nA_i-A_iX_n&=&\half (A_{n-i}A_i+A_iA_{n-i})A_i-\half A_i(A_{n-i}A_i+A_iA_{n-i})\\
&=&\half (A_{n-i}A_i-A_iA_{n-i})A_i+\half A_i (A_{n-i}A_i-A_iA_{n-i})\\
&=&\half (n-2i)(\E_nA_i+A_i\E_n).
\end{eqnarray*}
\end{ex}

\subsection{Representations}
Let us consider a well-known particular class of representations of $\K(1)$ called the
\textit{density representations}.
We will determine under what condition a given density representation is a representation
of $\A\K(1)$. 

The density representations are denoted by ${\mathcal F}_\lambda$,
where $\lambda\in\C$ is the parameter.
The basis in ${\mathcal F}_\lambda$ is $\{f_n,\,n\in\Z,\quad\phi_i,\,i\in\Z+\frac{1}{2}\}$
and the action of  ${\mathcal K}(1)$ is given by
$$
\begin{array}{rcl}
\chi_{x_n}(f_m) &=&
\left(m+\lambda{}n\right)f_{n+m},\\[8pt]
\chi_{x_n}(\phi_i) &=&
\left(i+(\lambda+\frac{1}{2})n\right)\phi_{n+i},\\[8pt]
\chi_{a_i}(f_n) &=&
\left(\frac{n}{2}+\lambda{}i\right)\phi_{i+n},\\[8pt]
\chi_{a_i}(\phi_j) &=&
2f_{i+j}.
\end{array}
$$

\begin{rem}

The adjoint representation of ${\mathcal K}(1)$ is
precisely the module ${\mathcal F}_{-1}$.

\end{rem}

\begin{prop}
\label{TDAdProf}
The ${\mathcal K}(1)$-module ${\mathcal F}_\lambda$ is a representation of ${\mathcal AK}(1)$
if and only if $\lambda=0$, or $\lambda=\frac{1}{2}$.
\end{prop}

\begin{proof}
Suppose that ${\mathcal F}_\lambda$ is a representation of ${\mathcal AK}(1)$.
By definition, the odd generators of ${\mathcal AK}(1)$ are represented by
$
\chi_{a_i}.
$
Surprisingly, one can check that operators on ${\mathcal F}_\lambda$ of the form:
$$
\textstyle
\frac{1}{j-i}\left(\chi_{a_i}\chi_{a_j}-\chi_{a_j}\chi_{a_i}\right)
$$
only depends on $i+j$ and not on the couple $(i,j)$.

We define the even generators by:
$$
\textstyle
\chi_{\varepsilon_n}:=
\frac{1}{j-i}\left(\chi_{a_i}\chi_{a_j}-\chi_{a_j}\chi_{a_i}\right)
$$
with $i+j=n$.
One then obtains
$$
\textstyle
\chi_{\varepsilon_n}(f_m) =
2\lambda\,f_{n+m},
\qquad
\chi_{\varepsilon_n}(\phi_i) =
(1-2\lambda)\phi_{n+i}.
$$
The relation
$$
\textstyle
\chi_{\varepsilon_n}\chi_{a_i}+\chi_{a_i}\chi_{\varepsilon_n}=\,\chi_{a_{n+i}}
$$
is always satisfied.
The operators $\chi_{\varepsilon_n}$ and $\chi_{\varepsilon_m}$ obviously commute, but
the relation
$$
\chi_{\varepsilon_n}\chi_{\varepsilon_m}=\chi_{\varepsilon_{n+m}}
$$
is true if and only if $\lambda=0,\frac{1}{2}$.
\end{proof}

\section{Appendix: the technical proofs}\label{ProofSect}

Theorem \ref{thmbrkt} and Proposition \ref{propbrkt} are formulated in \cite{Ovsienko}
without proofs.
These statements are crucial for the whole theory and their proofs are far of being evident.
We think that it is important to have these proofs in a written form.

\subsection{Proof of Proposition \ref{propbrkt}}
We will need to use all the axioms (LA1)-(LA3) of Lie antialgebras.

(i) The definition of $[a\odot b, c]$ is compatible with relations \eqref{relequiv}:

\noindent
The expression $a(bc)+b(ac)$ is symmetric in $a,b$. Therefore, one immediately has 
$$[a\odot b, c]=[b\odot a, c].$$
Now, we want to show that for $a,b,c \in \ga_1$, $x\in \ga_0$ we have
$$[a\odot  bx, c]=[ax\odot b, c].$$
That is, we want to show that the expression
$$(*)= \,a((bx)c)+(bx)(ac)-(ax)(bc)-b((ax)c)$$
cancels.

Using the identity (LA2) we rewrite the first and last terms:
\begin{eqnarray*}
a((bx)c)&=&a(\, x(bc)-b(xc)\,)\\[6pt]
b((ax)c)&=&b(\, x(ac)-a(xc)\,).
\end{eqnarray*}

So, now we have
$$(*)= a(x(bc))-a(b(xc))+(bx)(ac)-(ax)(bc)-b( x(ac))+b(a(xc)).$$

Using the identity (LA3) for the odd elements $a,b$ and $(xc)$ we can change the 2nd and last term of (*) to:
$$-a(b(xc))+b(a(xc))=(xc)(ab).$$

Finally, we have
$$(*)=a(x(bc))+(bx)(ac)-(ax)(bc)-b( x(ac))+(xc)(ba).$$
We can factor out the element $x$ in each term. 
Indeed, using  (LA1) we get
\begin{eqnarray*}
 a(x(bc))&=&2(a(bc))x\\[6pt]
 (bx)(ac)&=&(b(ac))x\\[6pt]
 (ax)(bc)&=&(a(bc))x\\[6pt]
  b(x(ac))&=&2(b(ac))x\\[6pt]
  (xc)(ab)&=&((ab)c)x.
\end{eqnarray*}
Replacing the above expressions in (*) we obtain
$$(*)=(a(bc))x-(b(ac))x+((ab)c)x.$$
And we can see that (*)=0 by the Jacobi identity (LA3).\\

(ii)  The definition of $[a\odot b, c\odot d]$ is compatible with relations \eqref{relequiv}:

\noindent
To see that the expression in \eqref{superbracket} is well-defined one expands it using some relations.
Using the Jacobi identity for the elements of $\ga_1$ as well as the relations \eqref{relequiv} in $\gog_\ga$ one can express $[a\odot b, c\odot d]$ in different ways.

\begin{lem}
One has the following equivalent expressions
\begin{equation}\label{supbkt2}
\begin{array}{rcl}
[a\odot b, c\odot d]&=&a(bc)\odot d + b(ad)\odot c + b(ac)\odot d +a(bd)\odot c \\[6pt]
&=& d(bc)\odot a + c(ad)\odot b + d(ac)\odot b + c(bd)\odot a .
\end{array}
\end{equation}
\end{lem}

\begin{proof}
To get the first expression of \eqref{supbkt2} one uses the Jacobi identity on the odd elements of $\ga_1$:
\begin{eqnarray*}
 b(ac)\odot d +a(bd)\odot c &=& \big(a(bc)+c(ab)\big)\odot d + \big(b(ad)-d(ab)\big) \odot c\\
 &=& a(bc)\odot d + b(ad)\odot c + c(ab)\odot d -d(ab)\odot c.
\end{eqnarray*}
The last two terms cancel by the relations \eqref{relequiv} in $\gog_\ga$. Thus,
\begin{eqnarray*}
a(bc)\odot d + b(ad)\odot c + b(ac)\odot d +a(bd)\odot c &=& 2\,a(bc)\odot d +2\, b(ad)\odot c.\end{eqnarray*}
This establishes the first identity  of \eqref{supbkt2}. The second expression is immediately deduced using the equivalence relations  \eqref{relequiv} in $\gog_\ga$.
\end{proof}

Now, from \eqref{supbkt2} it is easy to check that the expression of $[a\odot b, c\odot d]$ is symmetric in $a,b$ and  symmetric in $c,d$. So, the definition is compatible with the relations $a\odot b = b \odot a$ and $c\odot d = d \odot c$.

We also check the compatibility with the relations $ax\odot b = a \odot bx$ and $cx\odot d = c\odot dx$. Indeed, from \eqref{supbkt2} one can write
\begin{equation}\label{supbkt4}
\begin{array}{rcl}
[a\odot b, c\odot d]&=& [a\odot b,c]\odot d +  [a\odot b,d]\odot c.
\end{array}
\end{equation}
Thus, using point (i), we deduce the compatibility with $ax\odot b = a \odot bx$. 

Now using  the expressions of \eqref{supbkt2}  we can see that $[a\odot b, c\odot d]$ is skew-symmetric in $(a,b)$, $(c,d)$. In other words, one has
\begin{equation}\label{supbkt5}
\begin{array}{rcl}
[a\odot b, c\odot d]&=& -[c\odot d, a\odot b].
\end{array}
\end{equation}
Therefore, we also deduce the compatibility with $cx\odot d = c\odot dx$. 

We have proved Proposition \ref{propbrkt}.

\subsection{Proof of Theorem \ref{thmbrkt}}(i) One needs to check the property of antisymmety of the bracket:
$$[X,Y]=-(-1)^{\bar X \bar Y}[Y,X].$$
In the case of two odd elements or one odd and one even, this property is immediate from the definitions 
\eqref{superbracket}. In the case of two even elements this property has already been observed in \eqref{supbkt5}.\\

(ii) One checks the generalized Jacobi identity of the bracket
$$(-1)^{\bar X \bar Z}[[X,Y],Z]+(-1)^{\bar Y \bar X}[[Y,Z],X]+(-1)^{\bar Z \bar Y}[[Z,X],Y]=0.$$

- The case of three odd elements is immediate from the definitions \eqref{superbracket}.

- The case of two odd elements and one even is quite immediate. Indeed, we can rewrite \eqref{supbkt4} as 
$$
[a\odot b, [c, d]]= [\,[a\odot b,c], d]+  [\,[a\odot b,d], c],
$$
that is equivalent to Jacobi identity.

- The case of two even elements and one odd is less involved.

Denote by $J$ the expression
$$J=\big[\,[a\odot b, c\odot d], \,e\big]-\big[\, a\odot b,\,[ c\odot d,\,e]\big]+\big[\, c\odot d, \,[\,a\odot b,\,e]\big].$$
Using the first expression of \eqref{supbkt2} for $[a\odot b, c\odot d]$ and the definition \eqref{superbracket}, we can expand $J$ to
\begin{eqnarray*}
J&=&(a(bc))(de)+d((a(bc))e)+(b(ad))(ce)+c((b(ad))e)+\\
&&\quad(b(ac))(de)+d((b(ac))e)+(a(bd))(ce)+c((a(bd))e)\\[4pt]
&&-a(b(c(de)))-b(a(c(de)))-a(b(d(ce)))-b(a(d(ce)))\\[4pt]
&&+c(d(a(be)))+d(c(a(be)))+c(d(b(ae)))+d(c(b(ae))).\\
\end{eqnarray*}
One can split $J$ into symmetric expressions in $c,d$
$$J=J_1(a,b,c,d,e)+J_1(a,b,d,c,e)-J_2(a,b,c,d,e)-J_2(a,b,d,c,e),$$
where
\begin{equation*}
\begin{array}{rcl}
J_1(a,b,c,d,e)&=&(a(bc))(de)-a(b(c(de)))+(b(ac))(de)-b(a(c(de)))\\[6pt]
J_2(a,b,c,d,e)&=&c(e(b(ad)))-c(d(b(ae)))+c(e(a(bd)))-c(d(a(be))).\\
\end{array}
\end{equation*}
We show that $J_2(a,b,c,d,e)-J_1(a,b,c,d,e)=0$. By symmetry in $c,d$, this will imply that $J=0$.\\

All the terms involved are elements of the Lie antialgebra $\ga$. Once again we use all the axioms (LA1)-(LA3) to prove that the terms vanish.
One has
\begin{eqnarray*}
2(a(bc))(de)&=&a\big((bc)(de)\big)\\
&=&a(b(c(de)))+a((b(de))c)\\
&=&a(b(c(de)))+c((b(de))a)+(b(de))(ac)\\
&=&a(b(c(de)))-c(a(b(de)))+(b(ac))(de)\\
&=&a(b(c(de)))-c(a(b(de)))+(a(bc))(de)+(c(ab))(de).
\end{eqnarray*}
Thus, one has
\begin{equation*}
(a(bc))(de)-a(b(c(de)))=-c(a(b(de)))+(c(ab))(de).
\end{equation*}
By inverting the role of $a$ and $b$ in the above identity we also deduce
\begin{equation*}
(b(ac))(de)-b(a(c(de)))=-c(b(a(de)))+(c(ba))(de).
\end{equation*}
Hence, adding the two above identities we deduce
$$J_1(a,b,c,d,e)=-c(a(b(de)))-c(b(a(de))).$$
We have factored out $c$ in $J_1$. 
We then have
\begin{eqnarray*}
\dfrac{1}{c}(J_2-J_1)&=&e(b(ad))-d(b(ae))+e(a(bd))-d(a(be))+a(b(de))+b(a(de))\\
&=&e(b(ad))-\{(db)(ae)-(d(ae))b\}+\{ (ea)(bd)-(e(bd))a \}\\
&&-d(a(be))+a(b(de))+b(a(de))\\
&=&e(b(ad))+(d(ae))b-(e(bd))a-d(a(be))+a(b(de))-(a(de))b\\
&=&e(b(ad))+\big(d(ae)-a(de)\big)b+a\big(e(bd)+b(de)\big)-d(a(be))\\
&=&e(b(ad))+(e(ad))b+a(d(be))+(a(be))d\\
&=&(eb)(ad)+(ad)(be)\\
&=&0.
\end{eqnarray*}

The case of three even elements follows now from the fact that even elements are products of two odd elements and from the fact the Jacobi identities in the other cases hold. Indeed, 
denote by $K$
$$K=\big[\,[a\odot b, c\odot d], \,e\odot f\big]-\big[\, a\odot b,\,[ c\odot d,\,e\odot f]\big]+\big[\, c\odot d, \,[\,a\odot b,\,e\odot f]\big].$$

\noindent
We rewrite each terms using the property \eqref{supbkt4}.
\begin{eqnarray*}
\big[\,[a\odot b, c\odot d], \,e\odot f\big]&=& \big[\,[a\odot b, c\odot d], \,e\big] \odot f +\big[\,[a\odot b, c\odot d], \,f\big] \odot e
\end{eqnarray*}

\begin{eqnarray*}
\big[\, a\odot b,\,[ c\odot d,\,e\odot f]\big]&=&\big[\, a\odot b,\,[ c\odot d,\,e]\odot f+[ c\odot d,\,f]\odot e\big]\\
&=&\big[\, a\odot b,\,[ c\odot d,\,e]\big]\odot f+\big[\, a\odot b,f\big]\odot [ c\odot d,\,e]\\
&&+\big[\, a\odot b,\,[ c\odot d,\,f]\big]\odot e+\big[\, a\odot b,e\big]\odot [ c\odot d,\,f]
\end{eqnarray*}

\begin{eqnarray*}
\big[\, c\odot d,\,[ a\odot b,\,e\odot f]\big]&=&\big[\, c\odot d,\,[ a\odot b,\,e]\odot f+[ a\odot b,\,f]\odot e\big]\\
&=&\big[\, c\odot d,\,[ a\odot b,\,e]\big]\odot f+\big[\, c\odot d,f\big]\odot [ a\odot b,\,e]\\
&&+\big[\, c\odot d,\,[ a\odot b,\,f]\big]\odot e+\big[\, c\odot d,e\big]\odot [a\odot b,\,f]
\end{eqnarray*}
We replace in the expression of $K$, we get
\begin{eqnarray*}
K&=& \big[\,[a\odot b, c\odot d], \,e\big] \odot f-\big[\, a\odot b,\,[ c\odot d,\,e]\big]\odot f+\big[\, c\odot d,\,[ a\odot b,\,e]\big]\odot f\\
&&+ \big[\,[a\odot b, c\odot d], \,f\big] \odot e-\big[\, a\odot b,\,[ c\odot d,\,f]\big]\odot e+\big[\, c\odot d,\,[ a\odot b,\,f]\big]\odot e.
\end{eqnarray*}
From the previous case, we deduce that $K=0$.

Theorem \ref{thmbrkt} is proved.




\begin{thebibliography}{99}



\bibitem{Bavula}
V. Bavula, F. van Oystaeyen,
\textit{ The simple modules of the Lie superalgebra ${\rm osp}(1,2)$}.
J. Pure Appl. Algebra  150  (2000),  no. 1, 41--52. 

\bibitem{BG}
A. Braverman, D. Gaitsgory,
\textit{Poincar\'e-Birkhoff-Witt theorem for quadratic algebras of Koszul type}.  
J. Algebra  181  (1996),  no. 2, 315--328.

\bibitem{CNS}
L.Corwin, Y.Ne'eman, S.Sternberg,
\textit{Graded Lie algebras in mathematics and physics (Bose-Fermi symmetry)}
 Rev. Modern Phys.  47  (1975), 573--603


\bibitem{Kaplansky}
I. Kaplansky, unpublished preprints, available at www.justpasha.org/math/links/subj/lie/kaplansky/

\bibitem{Kac}
V. Kac,
\textit{Classification of simple $Z$-graded Lie superalgebras and simple Jordan superalgebras}.
Comm. Algebra  5  (1977), no. 13, 1375--1400. 

\bibitem{Martinez}
C. Martinez,
\textit{Simplicity of Jordan Superalgebras and Relations 
with Lie Structures }.
Irish Math. Soc. Bulletin 50 (2003), 97Ð116 .


\bibitem{MZ1}
C. Martinez, E. Zelmanov,
{\it Specializations of Jordan superalgebras},
Canad. Math. Bull. Vol. 45 (4), 2002 pp. 653Ð671 .

\bibitem{MZ2}
C. Martinez, E. Zelmanov,
{\it Unital bimodules over the simple Jordan superalgebra $D(t)$},
Trans. Amer. Math. Soc.  358  (2006),  no. 8, 3637--3649 (electronic).

\bibitem{McC}
K. McCrimmon,
{\it Kaplansky Superalgebras},
J. Algebra 164 (1994), 656Ð694. 

\bibitem{MG}
S. Morier-Genoud,
{\it Representations of $\asl_2$},
Intern. Math. Res. Notices., 2009.


\bibitem{Ovsienko}
V. Ovsienko,
{\it Lie antialgebras: pr\'emices},
arXiv:0705.1629.

\bibitem{T1}
M. N. Trushina,
{\it Irreducible representations of a certain Jordan superalgebra},
J. Algebra Appl.  4  (2005),  no. 1, 1--14. 


\end{thebibliography}
\end{document}